\documentclass[10pt]{amsart}
\usepackage[psamsfonts]{amssymb}
\usepackage{amsthm,amsmath}

\title{Associative polynomial functions over bounded distributive lattices}

\author{Miguel Couceiro}
\address{Mathematics Research Unit, FSTC, University of Luxembourg, 6, rue Coudenhove-Kalergi, L-1359 Luxembourg-Kirchberg, Luxembourg.} \email{miguel.couceiro[at]uni.lu}

\author{Jean-Luc Marichal}
\address{Mathematics Research Unit, FSTC, University of Luxembourg, 6, rue Coudenhove-Kalergi, L-1359 Luxembourg-Kirchberg, Luxembourg.} \email{jean-luc.marichal[at]uni.lu }

\date{February 17, 2010}

\begin{document}
\maketitle

\theoremstyle{plain}
\newtheorem{theorem}{Theorem}
\newtheorem{lemma}[theorem]{Lemma}
\newtheorem{proposition}[theorem]{Proposition}
\newtheorem{corollary}[theorem]{Corollary}
\newtheorem{fact}[theorem]{Fact}
\newtheorem{claim}{Claim}

\theoremstyle{definition}
\newtheorem{definition}[theorem]{Definition}
\newtheorem{example}[theorem]{Example}

\theoremstyle{remark}
\newtheorem{conjecture}{Conjecture}
\newtheorem{remark}{Remark}

\newcommand{\card}[1]{\ensuremath{\lvert{#1}\rvert}}
\newcommand{\vect}[1]{\ensuremath{\mathbf{#1}}} 
\newcommand{\co}[1]{\ensuremath{\overline{#1}}}
\def\median{\mathop{\rm median}\nolimits}
\newcommand{\N}{\mathbb{N}}
\newcommand{\R}{\mathbb{R}}
\newcommand{\ran}{\mathrm{ran}}
\def\med{\mathop{\rm med}\nolimits}

\begin{abstract}
The associativity property, usually defined for binary functions, can be generalized to functions of a given fixed arity $n\geqslant 1$ as well as to functions of
multiple arities. In this paper, we investigate these two generalizations in the case of polynomial functions over bounded distributive lattices and
present explicit descriptions of the corresponding associative functions. We also show that, in this case, both generalizations of associativity are essentially the same.
\end{abstract}


\noindent{\bf MSC classes:} 28B15, 39B72 (Primary) 06D05 (Secondary)

\section{Introduction}

Let $X$ be an arbitrary nonempty set. Throughout this paper, we regard vectors $\vect{x}$ in $X^n$ as $n$-strings over $X$.
The $0$-string or \emph{empty} string is denoted by
$\varepsilon$ so that $X^0=\{\varepsilon\}$. We denote by $X^*$ the set of all strings over $X$, that is, $X^*=\bigcup_{n\in \mathbb{N}}X^n$.
Moreover, we consider  $X^*$ endowed with concatenation for which we adopt the juxtaposition notation.  For instance, if $\vect{x}\in X^n$, $y\in X$, and $\vect{z}\in X^m$, then $\vect{x}y\vect{z}\in X^{n+1+m}$. Furthermore, for $\vect{x} \in X^m$,  we use the short-hand notation
$\vect{x}^n=\vect{x}\cdots \vect{x}\in X^{n\times m}$.
In the sequel, we will be interested both in functions of a given fixed arity (i.e., functions $f\colon X^n\to X$) as well as in functions defined on $X^*$, that is, of the form
$g\colon X^*\to X$.
Given a function $g\colon X^*\to X$, we denote by $g_n$ the restriction of $g$ to $X^n$, i.e.\ $g_n:=g|_{X^n}$. In this way, each function $g\colon X^*\to X$
can be regarded as a family $(g_n)_{n\in \mathbb{N}}$ of functions $g_n\colon X^n\to X$. We convey that $g_0$ is defined by $g_0(\varepsilon)=\varepsilon$.

In this paper, we are interested in the associativity property, traditionally considered on binary functions.
Recall that a function $f\colon X^2\to X$ is said to be \emph{associative} if $f(f(xy)z)=f(xf(yz))$ for every $x,y,z\in X$.
The importance of this notion is made clear by its natural interpretation.
Essentially, it expresses the fact that the order in which variables are bracketed is not relevant.
This algebraic property was extended to functions $f\colon X^n\to X$, $n\geqslant 1$, as well as to functions $g\colon X^*\to X$
in somewhat different ways.

A function $f\colon X^n\to X$ is said to be \emph{associative} if, for every $\vect{x}\vect{z},\vect{x}'\vect{z}'\in
X^{n-1}$ and every $\vect{y},\vect{y}'\in X^n$ such that $\vect{x}\vect{y}\vect{z}=\vect{x}'\vect{y}'\vect{z}'$, we have
$f(\vect{x}f(\vect{y})\vect{z})=f(\vect{x}'f(\vect{y}')\vect{z}')$. This generalization of associativity to $n$-ary functions goes back to D\"ornte \cite{Do1928} and led to the generalization of groups to $n$-groups (polyadic groups).\footnote{The first extensive study on polyadic groups was due to Post \cite{Post2}. This study was followed by several contributions towards the classification and description of $n$-groups and similar ``super-associative" structures; to mention a few, see \cite{Dudek95,Dudek01,DGG77,G,GG67,H63,MonkSioson71,Zup67}.} In a somewhat different context, this notion has been recently used to completely classify closed intervals made of equational classes of Boolean functions; see \cite{Cou08}.

On a different setting, associativity can be generalized to functions on $X^*$ as follows. We say that a function $g\colon X^*\to X$ is \emph{associative} if, for every
$\vect{x}\vect{y}\vect{z},\vect{x}'\vect{y}'\vect{z}'\in X^*$ such that $\vect{x}\vect{y}\vect{z}=\vect{x}'\vect{y}'\vect{z}'$, we have
$g(\vect{x}g(\vect{y})\vect{z})=g(\vect{x}'g(\vect{y}')\vect{z}')$. Alternative formulations of this definition appeared in the theory of aggregation functions, where the arity is not always fixed; see for instance \cite{BelPraCal07,GraMarMesPap09,KleMesPap00,Mar98}.

In general, the latter definition is more restrictive on the components $g_n$ of $g\colon X^*\to X$. For instance, the ternary real function $f(xyz)=x-y+z$ is associative but cannot be the ternary component of an associative function $g\colon \R^*\to \R$. Indeed, the equations
$$
g_2(g_2(xy)z)=g_2(xg_2(yz))=x-y+z
$$
have no solution,
for otherwise we would have $y=g_2(g_2(y0)0)$ and hence
$$
g_2(xy)=g_2(xg_2(g_2(y0)0))=g_2(g_2(xg_2(y0))0)=x-g_2(y0),
$$
which would imply $g_2(xy)=x-y$, a contradiction.

In this paper we show that, in the case of lattice polynomial functions, the two notions of associativity are essentially the same. More precisely, given a bounded distributive lattice $L$, we have that a polynomial function $f\colon L^n\to L$ is associative if and only if it is the $n$-ary component of some associative function $g\colon L^*\to L$; see Corollary~\ref{cor:main}. From this result and a characterization of polynomial functions given in \cite{CouMar0}, we derive a description of associative and range-idempotent polynomial functions $g\colon L^*\to L$ in terms of necessary and sufficient conditions. To this extent, in Section 2 we provide some preliminary results, which are then used in Section 3 to obtain explicit descriptions of those associative polynomial functions; see Theorems~\ref{thm:MainThm2}, \ref{thm:MainThm3}, and \ref{thm:CharLP23}.

\section{Preliminary results}

The following proposition provides useful reformulations of associativity of functions $g\colon X^*\to X$.

\begin{proposition}\label{prop:AltFormAs}
Let $g\colon X^*\to X$ be a function. The following assertions are equivalent:
\begin{enumerate}
\item[(i)] $g$ is associative.

\item[(ii)] For every $\vect{x}\vect{y}\vect{z}\in X^*$, we have $g(\vect{x}g(\vect{y})\vect{z})=g(\vect{x}\vect{y}\vect{z})$.

\item[(iii)] For every $\vect{x}\vect{y}\in X^*$, we have $g(g(\vect{x})g(\vect{y}))=g(\vect{x}\vect{y})$.

\end{enumerate}
\end{proposition}

\begin{proof}
We prove $(i)\Rightarrow (ii)$ simply by considering $\vect{x}\vect{y}\vect{z}\in X^*$, $\vect{x}'=\vect{x}\vect{y}\vect{z}$ and
$\vect{y}'=\vect{z}'=\varepsilon$. Also, we have $(ii)\Rightarrow (i)$ and $(ii)\Rightarrow (iii)$ trivially. Finally, let us prove that
$(iii)\Rightarrow (ii)$. First observe that $g(g(\vect{x}))=g(\vect{x})$ for every $\vect{x}\in X^*$. Therefore, for every
$\vect{x}\vect{y}\vect{z}\in X^*$, we have
$g(\vect{x}g(\vect{y})\vect{z})=g(g(\vect{x}g(\vect{y}))g(\vect{z}))=g(g(g(\vect{x})g(\vect{y}))g(\vect{z}))
=g(g(\vect{x}\vect{y})g(\vect{z}))=g(\vect{x}\vect{y}\vect{z})$.
\end{proof}

\begin{remark}\label{rem:AltDefAsG}
\begin{enumerate}
\item[(i)] Associativity of functions $g\colon X^*\to X$ was defined in \cite{KleMesPap00} and \cite{Mar98} as in assertions $(iii)$ and $(ii)$ of Proposition~\ref{prop:AltFormAs}, respectively. For a recent reference, see \cite{GraMarMesPap09}.

\item[(ii)] As observed in \cite{BelPraCal07}, associative functions $g\colon X^*\to X$ are completely determined by their unary and binary components. Indeed, for every $n\in\N$, $n>2$, and every $x_1,\ldots,x_n\in X$, we have
$$
g(x_1\cdots x_n)=g_2(g_2(\cdots g_2(g_2(x_1x_2)x_3)\cdots) x_n).
$$
\end{enumerate}
\end{remark}

A function $f\colon X^n\to X$ is said to be \emph{idempotent} if $f(x^n)=x$ for every $x\in X$. It is said to be \emph{range-idempotent}
\cite{GraMarMesPap09} if $f(f(\vect{x})^n)=f(\vect{x})$ for every $\vect{x}\in X^n$. Similarly, we say that a function $g\colon X^*\to X$ is
\emph{range-idempotent} if $g(g(\vect{x})^n)=g(\vect{x})$ for every $\vect{x}\in X^*$ and every integer $n\geqslant 1$.

\begin{lemma}\label{lemma:AsRanId}
Let $g\colon X^*\to X$ be an associative function. Then $g$ is range-idempotent if and only if $g(\vect{x}^n)=g(\vect{x})$ for every
$\vect{x}\in X^*$ and every integer $n\geqslant 1$.
\end{lemma}

\begin{proof}
For the sufficiency, simply observe that, for every $\vect{x}\in X^*$ and every $n\geqslant 1$, we have
$g(g(\vect{x})^n)=g(g(\vect{x}))=g(\vect{x})$. For the necessity, by repeated applications of Proposition~\ref{prop:AltFormAs} $(ii)$, we observe that, for every $\vect{x}\in X^*$ and every $n\geqslant 1$, we have
$g(\vect{x}^n)=g(g(\vect{x})^n)=g(\vect{x})$.
\end{proof}

\begin{lemma}\label{lemma:AsRanIdStrId}
Let $g\colon X^*\to X$ be an associative and range-idempotent function. Then, for every $\vect{x}y\vect{z}\in X^*$, we have
$g(\vect{x}g(\vect{x}y\vect{z})\vect{z})=g(\vect{x}y\vect{z})$.
\end{lemma}

\begin{proof}
Let $\vect{x}y\vect{z}\in X^*$. Using associativity and Lemma~\ref{lemma:AsRanId}, we have
$g(\vect{x}g(\vect{x}y\vect{z})\vect{z})=g(g(\vect{x}^2)y g(\vect{z}^2))=g(g(\vect{x})y g(\vect{z}))=g(\vect{x}y\vect{z})$.
\end{proof}

\section{Associative polynomial functions}

Let $L$ be a bounded distributive lattice, with 0 and 1 as bottom and top elements. In this section, we focus on (lattice) polynomial functions
$f\colon L^n\to L$, that is, functions which can be obtained as combinations of projections and constant functions using the lattice operations
$\wedge$ and $\vee$. As it is well known, these coincide exactly with those functions representable in disjunctive normal form (DNF).

More precisely, for $I\subseteq [n]=\{1,\ldots ,n\}$, let $\vect{e}_I\in\{0,1\}^n$ be the characteristic vector of $I$ and let $\alpha_f \colon 2^{[n]}\rightarrow L$
be the function given by $\alpha_f(I)=f(\vect{e}_I)$. Then
\begin{equation}\label{eq:LP-DNF}
f(\vect{x})=\bigvee_{I\subseteq [n]}\big(\alpha_f(I)\wedge \bigwedge_{i\in I} x_i\big).
\end{equation}
Moreover, by considering the function $\alpha^{*}_f \colon 2^{[n]}\rightarrow L$ defined by
\[
\alpha^{*}_f(I) =
\begin{cases}
\alpha_f(I), & \text{if $\bigvee_{J\varsubsetneq I}\alpha_f(J) < \alpha_f(I)$,} \\
0, & \text{otherwise,}
\end{cases}
\]
we obtain the `minimal' DNF representation of $f$ by replacing $\alpha_f(I)$ with $\alpha^*_f(I)$ in (\ref{eq:LP-DNF}). For further background,
see \cite{CouMar0, CouMar1}.

The following proposition provides a characterization of the $n$-ary polynomial functions. Recall that the ternary median function is the polynomial function $\mathrm{med}(x,y,z)=(x\vee y)\wedge (x\vee z)\wedge (y\vee z)$.

\begin{proposition}[{\cite{Marc}}]\label{Theorem:MAR34} A function $f\colon L^{n}\rightarrow L$ is a polynomial function if and only if $f(\vect{x}y\vect{z})=\med\big(f(\vect{x}0\vect{z}),y, f(\vect{x}1\vect{z})\big)$ for every $\vect{x}y\vect{z}\in L^n$.
\end{proposition}

The following theorem restricts the disjunctive normal form of $n$-ary associative polynomial functions. We first consider a lemma which follows immediately from Proposition~\ref{Theorem:MAR34}.

\begin{lemma}\label{lemma:MedDec}
Let $f\colon L^n\to L$ be a polynomial function and let $I\varsubsetneq [n]$, $J\subseteq [n]$, and $k\in [n]\setminus I$. Then, for $\vect{x}0\vect{z}=\vect{e}_I$ with $\vect{x}\in L^{k-1}$, we have
$$
f(\vect{x}f(\vect{e}_J)\vect{z})=\mathrm{med}\big(\alpha_f(I),\alpha_f(J),\alpha_f(I\cup\{k\})\big).
$$
\end{lemma}

\begin{theorem}\label{thm:MainThm2}
Let $f\colon L^n\to L$ be a polynomial function. If $f$ is associative, then
\begin{equation}\label{eq:MainThm2}
f(\mathbf{x}) = a_n\vee(b_n\wedge x_1)\vee\Big(\bigvee_{i=1}^{n}(b_n\wedge c_n\wedge x_i)\Big)\vee (c_n\wedge
x_n)\vee\Big(d_n\wedge\bigwedge_{i=1}^nx_i\Big),
\end{equation}
where $a_n=f(0^n)$, $b_n=f(10^{n-1})$, $c_n=f(0^{n-1}1)$, and $d_n=f(1^n)$.
\end{theorem}

\begin{proof} Let $f\colon L^n\to L$ be an associative polynomial function. Without loss of generality, we may assume that $n\geqslant 3$ for if $n=1$ or $n=2$, then it is trivial.

First, we show that $\alpha^*_f(I)=0$ whenever $I\subset [n]$ and $1<|I|<n$. To reach a contradiction, suppose that there is $I\subset [n]$ with $1<|I|<n$ and such that $\alpha^*_f(I)\neq 0$. Let $j\in [n]$ be the least such that $j\notin I$.
\begin{itemize}
\item[(i)] Suppose $j=1$. Using Lemma~\ref{lemma:MedDec}, from $f(f(0^n)\vect{x})=f(0^{n-1}f(\vect{e}_I))$, where $0\vect{x}=\vect{e}_I$, we obtain $\alpha_f(I)=\alpha_f(\{n\})\wedge \alpha_f(I)$, that is, $\alpha_f(I) \leqslant \alpha_f(\{n\})$. This implies $n\notin I$, for otherwise we would have $\alpha_f(I)=\alpha_f^*(I)>\alpha_f(\{n\})$. Similarly, from  $f(\vect{y}f(0^n))=f(f(\vect{e}_I)0^{n-1})$, where $\vect{y}0=\vect{e}_I$, we obtain $\alpha_f(I)=\alpha_f(\{1\})\wedge\alpha_f(I)$. Finally, let $k\in [n]\setminus\{1\}$ be the least such that $k\in I$ and take $\vect{z}\in L^n$ such that $0^{k-1}\vect{z}0^{n-k}=\vect{e}_I0^{n-1}$. Then
    \begin{eqnarray*}
    f(f(\vect{e}_I)0^{n-1}) &=& \alpha_f(\{1\})\wedge\alpha_f(I)\\
    &=& \alpha_f(I)=\alpha^*_f(I)>\alpha_f(\{k\})\geqslant f(0^{k-1}f(\vect{z})0^{n-k}),
    \end{eqnarray*}
    which contradicts associativity.

\item[(ii)] Suppose $j>1$. Take $\vect{x}\in L^{n-j}$ such that $\vect{e}_I=1^{j-1}0\vect{x}$. On the one hand, we have $f(1^{j-1}f(0^{n})\vect{x})=\alpha_f(I)$. On the other hand, we obtain
    $$
    f(f(1^{j-1}0^{n-j+1})0^{j-1}\vect{x})
    \begin{cases}
    =\alpha_f(\{1\}), & \mbox{if $\vect{x}=0^{n-j}$,}\\
    <\alpha_f(I), & \mbox{otherwise},
    \end{cases}
    $$
    which contradicts associativity.
\end{itemize}

Next, we show that $\alpha_f(\{1\})\wedge \alpha_f(\{n\})\geqslant \alpha_f(\{i\})$ for every $i\in [n]\setminus \{1,n\}$.
Let $i\in [n]\setminus \{1,n\}$. We have $f(f({0}^{i-1}1{0}^{n-i}){0}^{n-1})=\alpha_f(\{1\})\wedge \alpha_f(\{i\})$ and $f({0}^{i-1}1{0}^{n-i-1}f({0}^{n}))=\alpha_f(\{i\})$. By associativity, it follows that $\alpha_f(\{i\})\leqslant \alpha_f(\{1\})$.
Similarly, we can verify that $\alpha_f(\{i\})\leqslant \alpha_f(\{n\})$, for every $i\in [n]\setminus \{1,n\}$.

Finally, we show that, for every $i\in [n]\setminus \{1,n\}$, we have $\alpha_f(\{1\})\wedge \alpha_f(\{n\})=\alpha_f(\{i\})$.
This will be enough to show that (\ref{eq:MainThm2}) holds, with $a=\alpha_f(\varnothing)$, $b=\alpha_f(\{1\})$, $c=\alpha_f(\{n\})$, and $d=\alpha_f([n])$.
So let $i\in [n]\setminus \{1,n\}$. On the one hand, we have
\begin{eqnarray*}
\lefteqn{f(f(0^{n-1}1)0^{i-2}10^{n-i})}\\ &=& \alpha_f(\{i\})\vee\big(\alpha_f(\{n\})\wedge\alpha_f(\{1,i\})\big)\\
&=& \alpha_f(\{n\})\wedge\alpha_f(\{1,i\})\qquad (\mbox{since $\alpha_f(\{i\})\leqslant\alpha_f(\{n\})$})\\
&=& \alpha_f(\{n\})\wedge\big(\alpha_f(\{1\})\vee\alpha_f(\{i\})\big)\qquad (\mbox{since $\alpha^*_f(\{1,i\})=0$})\\
&=& \alpha_f(\{1\})\wedge\alpha_f(\{n\})\qquad (\mbox{since $\alpha_f(\{i\})\leqslant\alpha_f(\{1\})$}).
\end{eqnarray*}
On the other hand, we have
\begin{eqnarray*}
f(0^{i-1}f(0^{n-i}10^{i-2}1)0^{n-i})&=& \alpha_f(\{i\})\wedge\alpha_f(\{n-i+1,n\})\\
&=& \alpha_f(\{i\})\wedge\alpha_f(\{n\})=\alpha_f(\{i\}),
\end{eqnarray*}
and the proof is now complete.
\end{proof}

\begin{remark}
\begin{enumerate}
\item[(i)] We observe that equation (\ref{eq:MainThm2}) can be rewritten in a more symmetric way as
$$
f(\mathbf{x}) = \mathrm{med}\Big(a_n,(b_n\wedge x_1)\vee\mathrm{med}\Big(\bigwedge_{i=1}^nx_i,b_n\wedge c_n,\bigvee_{i=1}^nx_i\Big)\vee (c_n\wedge
x_n),d_n\Big).
$$
This formula reduces to $f(\mathbf{x}) = \mathrm{med}\big(a_n,\mathrm{med}\big(\bigwedge_{i=1}^nx_i,b_n,\bigvee_{i=1}^nx_i\big),d_n\big)$ as soon as $f$ is a symmetric function (i.e., invariant under permutation of its variables).
\item[(ii)] A \emph{term function} $f\colon L^n\to L$ is a polynomial function satisfying $\alpha_f(I)\in\{0,1\}$ for every $I\subseteq [n]$. By Theorem~\ref{thm:MainThm2}, the only associative term functions $f\colon L^n\to L$ are $\vect{x}\mapsto x_1$, $\vect{x}\mapsto x_n$, $\vect{x}\mapsto\bigwedge_{i=1}^nx_i$, and $\vect{x}\mapsto\bigvee_{i=1}^nx_i$.
\end{enumerate}
\end{remark}

We say that a function $g\colon L^*\to L$ is a \emph{polynomial function} if every $g_n$, $n\geqslant 1$, is a polynomial function. The
following theorem yields a description of associative polynomial functions $g\colon L^*\to L$.

\begin{theorem}\label{thm:MainThm3}
A polynomial function $g\colon L^*\to L$ is associative if and only if $g_1(x) = a_1\vee (d_1\wedge x)$ and, for $n\geqslant 2$,
$$
g_n(\mathbf{x}) = a_2\vee(b_2\wedge x_1)\vee\Big(\bigvee_{i=1}^{n}(b_2\wedge c_2\wedge x_i)\Big)\vee (c_2\wedge
x_n)\vee\Big(d_2\wedge\bigwedge_{i=1}^nx_i\Big),
$$
where $a_1=g_1(0)$, $d_1=g_1(1)$, $a_2=g_2(0^2)$, $b_2=g_2(10)$, $c_2=g_2(01)$, $d_2=g_2(1^2)$, $a_1\leqslant a_2$, and $d_2\leqslant d_1$.
\end{theorem}

\begin{proof}
The sufficiency can be easily verified using Proposition~\ref{prop:AltFormAs}.

Let us establish the necessity. Since every function $g_n$, $n\geqslant 1$, is an
associative polynomial function, by Theorem~\ref{thm:MainThm2} this function has the form given in the right-hand side of (\ref{eq:MainThm2}), with $a_n=g_n(0^n)$,
$b_n=g_n(10^{n-1})$, $c_n=g_n(0^{n-1}1)$, and $d_n=g_n(1^n)$. By associativity and Proposition~\ref{Theorem:MAR34}, for every $n\geqslant 3$, we have
$$
g_n(10^{n-1})=g_2(g_{n-1}(10^{n-2})0)=\mathrm{med}\big(g_2(0^2),g_{n-1}(10^{n-2}),g_2(10)\big).
$$
By reasoning recursively, it can be verified that $b_n=b_2$ for every $n\geqslant 3$. Similarly, it can be shown that $a_n=a_2$, $c_n=c_2$, and $d_n=d_2$ every $n\geqslant 3$. Finally, By Propositions~\ref{prop:AltFormAs} and \ref{Theorem:MAR34}, we have
$$
g_2(x^2)=g_1(g_2(x^2))=\mathrm{med}\big(g_1(0),g_2(x^2),g_1(1)\big),
$$
which shows that $a_1\leqslant a_2$ and $d_2\leqslant d_1$.
\end{proof}

Even though associativity for functions $g\colon
L^*\to L$ seems more restrictive on their components $g_n$ than associativity for functions of a given fixed arity, from Theorems~\ref{thm:MainThm2}
and \ref{thm:MainThm3} it follows that associativity for polynomial functions $f\colon L^n\to L$ naturally extends componentwise to polynomial functions $g\colon L^*\to L$.

\begin{corollary}\label{cor:main}
Let $f\colon L^n\to L$ be a polynomial function. Then $f$ is associative if and only if there is an associative polynomial function $g\colon
L^*\to L$ such that $g_n = f$.
\end{corollary}

We now provide a characterization of the associative and range-idempotent polynomial functions $g\colon L^*\to L$ in terms of necessary and sufficient conditions. To this extent, we present a characterization of the $n$-ary polynomial functions given in \cite{CouMar0}. Recall that $S\subseteq L$ is \emph{convex} if, for every $y\in L$, we have that $x\leqslant y\leqslant z$ implies $y\in S$ whenever $x,z\in S$.

\begin{proposition}\label{thm:CharLP}
A function $f\colon L^n\to L$ is a polynomial function if and only if it satisfies
\begin{enumerate}
\item[(a)] $f$ is nondecreasing,

\item[(b)] for every $\vect{xz}\in L^{n-1}$, the function $y\mapsto f(\vect{x}y\vect{z})$ preserves $\wedge$ and $\vee$,

\item[(c)] $\{f(\vect{x}y\vect{z}):y\in L\}$ is convex for every $\vect{xz}\in L^{n-1}$,

\item[(d)] $\{f(\vect{x}):\vect{x}\in L^n\}$ is convex,

\item[(e)] $f(\vect{x}f(\vect{x}y\vect{z})\vect{z})=f(\vect{x}y\vect{z})$ for every $\vect{x}y\vect{z}\in L^n$.
\end{enumerate}
\end{proposition}

\begin{theorem}\label{thm:CharLP23}
Let $g\colon L^*\to L$ be an associative function. The following assertions are equivalent:
\begin{enumerate}
\item[(i)] $g$ is a range-idempotent polynomial function.

\item[(ii)] $g$ is a polynomial function satisfying $g_2(0^2)=g_1(0)$ and $g_2(1^2)=g_1(1)$.

\item[(iii)] $g$ is range-idempotent and, for every $n\geqslant 1$, the function $g_n$ satisfies
\begin{enumerate}
\item[(a)] $g_n$ is nondecreasing,

\item[(b)] for every $\vect{xz}\in L^{n-1}$, the function $y\mapsto g_n(\vect{x}y\vect{z})$ preserves $\wedge$ and $\vee$,

\item[(c)] $\{g_n(\vect{x}y\vect{z}):y\in L\}$ is convex for every $\vect{xz}\in L^{n-1}$,

\item[(d)] $\{g_n(\vect{x}):\vect{x}\in L^n\}$ is convex.
\end{enumerate}
\end{enumerate}
\end{theorem}

\begin{proof}
The equivalence $(i)\Leftrightarrow (iii)$  follows from Lemma~\ref{lemma:AsRanIdStrId} and Proposition~\ref{thm:CharLP}. To see that $(i)\Rightarrow (ii)$, we use Lemma~\ref{lemma:AsRanId} and Theorem~\ref{thm:MainThm3}. Finally, to show that $(ii)\Rightarrow (i)$, by Lemma~\ref{lemma:AsRanId} and Theorem~\ref{thm:MainThm3}, we only need to show that $g$ satisfies $g(\vect{x}^n)=g(\vect{x})$ for every $\vect{x}\in X^*$, which is immediate.
\end{proof}

In the case when $L$ is a bounded chain, that is, a totally ordered bounded lattice, the condition $(b)$ in Theorem~\ref{thm:CharLP23} (iii) becomes redundant in the presence of condition $(a)$. Moreover, it was shown in \cite[Lemma~18]{CouMar2} that conditions $(a)$ and $(c)$ together imply condition $(d)$. We then have the following corollary.

\begin{corollary}\label{cor:CharLP23}
Assume that $L$ is a bounded chain and let $g\colon L^*\to L$ be an associative function. The following assertions are equivalent:
\begin{enumerate}
\item[(i)] $g$ is a range-idempotent polynomial function.

\item[(ii)] $g$ is a polynomial function satisfying $g_2(0^2)=g_1(0)$ and $g_2(1^2)=g_1(1)$.

\item[(iii)] $g$ is range-idempotent and, for every $n\geqslant 1$, the function $g_n$ is nondecreasing and $\{g_n(\vect{x}y\vect{z}):y\in L\}$ is convex for every $\vect{xz}\in L^{n-1}$.
\end{enumerate}
\end{corollary}

\begin{remark}
In the special case of real intervals, i.e., when $L = [a,b]$ for reals $a\leqslant b$, the convexity condition can be
replaced with continuity of $g_n$ in Corollary~\ref{cor:CharLP23} (iii).\footnote{Indeed, for nondecreasing functions $f\colon [a,b]^n\to [a,b]$, this convexity condition is equivalent to continuity of $f$ in each variable and hence to continuity of $f$.} The more general case when $L$ is a connected ordered topological space was considered by
Fodor \cite{For96} who obtained an explicit description of those nondecreasing binary functions which are idempotent, continuous, and associative.
\end{remark}

\end{document}